\documentclass[letterpaper, 10 pt, conference]{ieeeconf}  

\usepackage{amsmath, amssymb, bbm, xspace}
\usepackage{cite}

\def\diag{\mathop{\hbox{\rm diag}}}

\def\grad{\nabla}

\def\spose#1{\hbox to 0pt{#1\hss}}

\def\text #1{\hbox{\quad#1\quad}}

\def\nthinsp{\mskip -2   mu}

\def\F{_{\scriptscriptstyle F}}

\def\k{_k}

\def\N{\mathbb{N}}

\def\superstar{^{\raise 0.5pt\hbox{$\nthinsp *$}}}
\def\SUPERSTAR{^{\raise 0.5pt\hbox{$*$}}}

\def\lamstarT {\lambda^{\raise 0.5pt\hbox{$\nthinsp *$}T}}

\def\mubar{\skew3\bar \mu}

\def\sigmabar{\bar\sigma}

\def\hbar{\skew{4.2}\bar h}

\def\Qbar{\bar Q}

		\def\bkE{{\rm I\kern-.17em E}}
		\def\bk1{{\rm 1\kern-.17em l}}
		\def\bkD{{\rm I\kern-.17em D}}
		\def\bkR{{\rm I\kern-.17em R}}
		\def\bkP{{\rm I\kern-.17em P}}
		\def\bkY{{\bf \kern-.17em Y}}
		\def\bkZ{{\bf \kern-.17em Z}}

		\def\beq{\begin{eqnarray}}
		\def\bc{\begin{center}}
		\def\be{\begin{enumerate}}
		\def\bi{\begin{itemize}}
		\def\bS{\begin{slide}}
		\def\ec{\end{center}}
		\def\ee{\end{enumerate}}
		\def\ei{\end{itemize}}
		\def\eS{\end{slide}}
		\def\eeq{\end{eqnarray}}

	\def\cp2problem#1#2#3#4{\fbox
		 {\begin{tabular*}{0.9\textwidth}
			{@{}l@{\extracolsep{\fill}}l@{\extracolsep{6pt}}l@{\extracolsep{\fill}}c@{}}
				#1 & & $#4 $ 
			\end{tabular*}}}

		\renewcommand{\emph}[1]{\textbf{#1}}

		\def\bkE{{\rm I\kern-.17em E}}
		\def\bk1{{\rm 1\kern-.17em l}}
		\def\bkD{{\rm I\kern-.17em D}}
		\def\bkR{{\rm I\kern-.17em R}}
		\def\bkP{{\rm I\kern-.17em P}}
		
		\def\bkZ{{\bf{Z}}}

\newcommand {\beeq}[1]{\begin{equation}\label{#1}}
\newcommand {\eeeq}{\end{equation}}
\newcommand {\bea}{\begin{eqnarray}}
\newcommand {\eea}{\end{eqnarray}}

\def\texitem#1{\par\smallskip\noindent\hangindent 25pt
               \hbox to 25pt {\hss #1 ~}\ignorespaces}

\newtheorem{definition}{Definition}{\it}{}
{\it}{}
{\it}{}
{\it}{}
\newtheorem{lemma}{Lemma}{\it}{}
\newtheorem{theorem}{Theorem}{\it}{}
\newtheorem{remark}{Remark}{\it}{}
{\it}{}
\newtheorem{standing}{Standing Assumption}

\def\diag{{\rm diag}}

\def\R{\mathbb{R}}

\def\F{\mathcal{F}}

\def\argmin{\mathop{\rm argmin}}

\def\x{\bs{x}}
\def\xhb{{\hat{\bs{x}}}}

\def\y{\bs{y}}

\def\F{{\bs{F}}}

\def\ber{\begin{eqnarray}}
\def\eer{\end{eqnarray}}
\def\bers{\begin{eqnarray*}}
	\def\eers{\end{eqnarray*}}
\def\be{\begin{equation}}
\def\ee{\end{equation}}
\def\1{{\bf 1}}

\newcommand{\bs}{\boldsymbol}
\newcommand{\mc}{\mathcal}

\newcommand{\proj}{\mathrm{proj}}
\newcommand{\col}{\mathrm{col}}

\renewcommand{\emph}{\textit}

\newcommand{\0}{\bs 0}

\def\k{{k \in \N}}

\def\W{\bs{W}}
\def\V{\bs{V}}

\usepackage{amsmath} 
\usepackage{amsfonts,mathrsfs}
\usepackage{amssymb}
\usepackage{mathtools}
\usepackage{color}
\usepackage{graphicx}
\usepackage{epsfig}
\usepackage{algorithm}
\usepackage{algorithmicx}
\usepackage[acronym]{glossaries}
\usepackage{subfigure}
\usepackage{cancel}
\usepackage{empheq}
\usepackage{placeins}

\let\labelindent\relax
\usepackage{enumitem}

\makeglossaries
\newacronym{GNEP}{GNEP}{generalized Nash equilibrium problem}
\newacronym{NE}{NE}{Nash equilibrium}
\newacronym{NEP}{NEP}{Nash equilibrium problem}
\newacronym{GNE}{GNE}{generalized Nash equilibrium}
\newacronym{v-GNE}{v-GNE}{variational \gls{GNE}}
\newacronym{ISS}{ISS}{Input-to-state-stable}
\newacronym{PPPA}{PPPA}{preconditioned proximal-point algorithm}
\newacronym{PPA}{PPA}{proximal-point algorithm}
\newacronym{VI}{VI}{variational inequality}
\newacronym{GAE}{GAE}{generalized aggregative equilibrium}
\newacronym{v-GAE}{v-GAE}{variational \gls{GAE}}
\newacronym{KKT}{KKT}{Karush-Kuhn-Tucker}
\newacronym{FQNE}{FQNE}{firmly quasinonexpansive}
\newacronym{FNE}{FNE}{firmly nonexpansive}
\newacronym{PF}{PF}{Perron-Frobenius}
\DeclareSymbolFont{myletters}{OML}{ztmcm}{m}{it}
\DeclareMathSymbol{\uplambda}{\mathord}{myletters}{"15}

\def\QEDhereeqn{\eqno\let\eqno\relax\let\leqno\relax\let\veqno\relax\hbox{\QED}}
\def\QEDopenhereeqn{\eqno\let\eqno\relax\let\leqno\relax\let\veqno\relax\hbox{\QEDopen}}
\newcommand\numberthis{\addtocounter{equation}{1}\tag{\theequation}}
\def\Q{\bs{Q}}

\def\FQ{\bs{\bar{F}}}
\def\fq{\bar{F}}
\def\barmu{\bar{\mu}}
\def\barell{\bar{\ell}}

\def\tildeell{\tilde{\ell}}
\def\eigmax#1{\uplambda_{\textnormal{max}}(#1)}
\def\eigmin#1{\uplambda_{\textnormal{min}}(#1)}
\def\eig{\uplambda}
\def\sigmamax#1{\sigma_{\textnormal{max}}(#1)}

\def\H{\mc{H}}
\IEEEoverridecommandlockouts 
\title{\LARGE \bf
Nash equilibrium seeking under partial-decision information over directed communication networks
}

\author{Mattia Bianchi and Sergio Grammatico
\thanks{The authors are with the Delft Center for Systems and Control (DCSC), TU Delft, The Netherlands.
	E-mail addresses: \texttt{\{m.bianchi, s.grammatico\}@tudelft.nl}. This work was partially supported by NWO under research project OMEGA (grant n. 613.001.702) and by the ERC under research project COSMOS (802348).} 
}

\begin{document}

\maketitle
\thispagestyle{empty}
\pagestyle{empty}

\begin{abstract}
	We consider the Nash equilibrium problem in a partial-decision information scenario. Specifically, each agent can only receive information from some neighbors via a communication network, while its  cost function  depends on the strategies of possibly all agents. In particular, while the existing methods assume undirected or balanced communication, in this paper we allow for non-balanced, directed graphs.
	 We propose a fully-distributed pseudo-gradient scheme, which is guaranteed to converge  with linear rate to a Nash equilibrium, under strong monotonicity and Lipschitz continuity of the game mapping. Our algorithm requires global knowledge of the communication structure, namely of the Perron-Frobenius eigenvector of the adjacency matrix and of a certain constant related to the graph connectivity. Therefore,  we adapt the procedure to setups where the network is not known in advance, by computing the  eigenvector online and by means of vanishing step sizes. 
\end{abstract}

\section{Introduction}
Game theory is a powerful tool to model and control the decision-making process of  selfish
agents, that aim at optimizing their individual,
but inter-dependent, objective functions. 
This scenario arises
 in  several relevant engineering applications, such as
congestion control in traffic networks \cite{Barrera2015}, smart-grid management \cite{Saad2012}, 
demand
response in competitive markets \cite{Li_Chen_Dahleh_2015} and  analysis of social dynamics \cite{Ghaderi_2014}.
Often, the goal (either of the agents or of a coordinator 
that pursues network regulation by imposing incentives or behavioral rules)
is the attainment of a \gls{NE}, a joint strategy 
from which 
it is not convenient for any agent to unilaterally deviate.

In fact, a recent part of the literature focuses on designing distributed \gls{NE} seeking algorithms, where the computational effort is partitioned among the agents \cite{Shamma_Arslan_2005,DePersisGrammaticoTAC2020,BelgioiosoGrammatico_ECC_2018}.
Nonetheless, typically these methods still  assume the presence of a central coordinator that can broadcast some data -- for instance, the average of all the agents' strategies, in the case of aggregative games \cite{BelgioiosoGrammatico_ECC_2018}. Unfortunately, this requirement is impractical in some domains \cite{SwensonKarXavier_FictitiousPlay_2015}.
To overcome this limitation, we consider \emph{fully-distributed} schemes, where the agents only rely on the information locally exchanged
over a network, via peer-to-peer communication. In particular, the main challenge is that the cost function of each agent may depend on the strategies of some other  non-neighboring agents. One example is the Cournot competition model described in \cite{Koshal_Nedic_Shanbag_2016}, where the profit of each of a group of firms depends not only on its own production, but also on the total supply, a quantity not directly accessible by any of the firms. 
To remedy  the lack of knowledge, each agent can estimate and eventually reconstruct the strategies of all the competitors (or an aggregation value),  based on the data received from its neighbors. 

Such a \emph{partial-decision information} setup has only been introduced very recently. Most of the available results resort to (projected) pseudo-gradient and consensus dynamics \cite{Koshal_Nedic_Shanbag_2016,TatarenkoShiNedic_CDC2018,DePersisGrammatico2018,YeHu2017,BelgioiosoNedicGrammatico2020,Pavel2018}. Alternatively, schemes based on a proximal-point iteration were studied in \cite{Bianchi_arXiv_TAC20_PPP}; a fully-distributed fictitious play algorithm was proposed in \cite{SwensonKarXavier_FictitiousPlay_2015}.
These approaches assume undirected communication, which might be unrealistic, e.g., in wireless systems, if the agents send signals at different power levels, implying unilateral transmission capability. Fewer works deal with asymmetric networks.
 Under the assumption of balanced weights, continuous-time dynamics were proposed in \cite{DengNian2019} for aggregative games;
 most recently, we also addressed generally-coupled-cost  games via a fixed-step forward-backward method \cite{Bianchi_LCSS2020}.
 To the best of our knowledge, the only  discrete-time \gls{NE} seeking algorithm that takes into account non-balanced  digraphs is the asynchronous gossip-based scheme in \cite{SalehisadaghianiPavel2017_nondoublystochastic}.

Even in the context of distributed optimization, most algorithms are designed with doubly stochastic adjacency matrices, which enjoy several convenient properties, not least that the  average of the agents' estimates is preserved over time. However, doubly stochastic weights cannot be easily assigned over directed networks. An alternative is to rely on column stochastic graphs, which  maintain the average invariance and only require the agents to know their out-degree.
Yet, this  is impractical in setups where the agents broadcast some information, but ignoring which of the other nodes can receive it; or if some of the communication links can fail. In contrast, distributed design of row stochastic matrices is straightforward, as it suffices for each agent to locally assign appropriate weights to the incoming information. However, the use of  row stochastic graphs comes with technical challenges, since many properties of doubly stochastic matrices are lost.  Of major interest for this work is the approach in \cite{XiMai_DistributedOptimizationonRowStochastic}: to correct the imbalance caused by employing row stochastic weights, the algorithm exploits the information contained in the \gls{PF} eigenvector of the adjacency matrix, which is computed online.  

\emph{Contribution:} Motivated by the above, we design the first synchronous, fully-distributed algorithm to compute a  \gls{NE} over directed non-balanced communication networks. Our contributions are summarized as follows:
{
\begin{itemize}[leftmargin=*,topsep=0pt]
	\item We prove  that any row stochastic primitive matrix with positive diagonal enjoys a contractivity property, in a Hilbert space weighted by its \gls{PF} eigenvector. We later exploit this general result to prove convergence of our equilibrium seeking dynamics (§\ref{sec:mathbackground});
\item We design a fully-distributed, fixed-step gradient algorithm to seek a \gls{NE} over strongly connected 
directed graphs, which is guaranteed to converge with liner rate under strong monotonicity of the game mapping. In our method, the pseudo-gradient component is divided by the entries of the \gls{PF} eigenvector of the network. Although this technique has already been adopted in distributed optimization  \cite{XiMai_DistributedOptimizationonRowStochastic}, we give a new, powerful, monotone-operator-theoretic interpretation, which greatly simplifies our analysis (§\ref{sec:math:A});
	\item  We show that convergence is retained even if the graph is not known in advance and the \gls{PF} eigenvector is computed online, provided that a small-enough step size is chosen. Since computing the upper bound distributedly can be troublesome, we also provided convergence guarantees for  vanishing steps (§\ref{sec:math:B}).
\end{itemize}
}

\smallskip
\emph{Basic notation}: 
$\R_{>0}$ denotes the set of positive real numbers.
$\0_n$ ($\1_n$) $\in \R^n$  denotes the vector with all elements equal to $0$ ($1$); $I_n\in\R^{n\times n}$ denotes an identity matrix; we may omit the subscripts if there is no ambiguity. $e_i^n\in\R^n$ denotes  a vector with all elements equal to $0$ except the $i$-th element, which is $1$. For a  function $g:\R^n \rightarrow \R$, $\nabla_{\! \! x} g(x)$ denotes its gradient.  For a matrix $A \in \R^{m \times n}$, $[A]_{i,j}$ represents the element on row $i$ and column $j$; $\sigma_{\textnormal{min}}(A)=\sigma_1(A)\leq\dots\leq\sigma_n(A)=:\sigma_{\textnormal{max}}(A)$ are its singular values.
 If $A\in\R^{n\times n}$ is symmetric, $\uplambda_{\textnormal{min}}(A)=\uplambda_1(A)\leq\dots\leq\uplambda_n(A)=:\uplambda_{\textnormal{max}}(A)$ denote its eigenvalues; $A \succ 0$  stands for  positive definite matrix.
$\otimes$ denotes the Kronecker product.
$\diag(A_1,\dots,A_N)$ denotes the block diagonal matrix with $A_1,\dots,A_N$ on its diagonal. 
Given $P\succ 0$, $\langle x, y \rangle _P=x^\top P y$ and
$\|{x}\|_{P}=\sqrt{x^\top P x}$ denote the $P$-weighted Euclidean inner product and norm, respectively;
 $\|A\|_P:=\sup_{x\neq \0} \! \textstyle \frac{\|Ax\|_P}{\|x\|_P}$ is the $P$-induced norm of $A\in \R^{n\times n}$; 
 we omit the subscripts if $P=I$. $\H_P:=(\R^n,\langle \cdot \, , \cdot \rangle_P)$ is the  Hilbert space obtained by endowing $\R^n$ with the $P$-weighted inner product.

\emph{Operator-theoretic notation}: 
 An operator $\mathcal{A} : \R^n \rightarrow \R^n$ is 
($\mu$-strongly) monotone in $\H_P$ if $\langle \mc{A}(x)-\mc{A}(y) ,  x-y\rangle _P \geq 0 \, (\geq \mu \|x-y\|_P^2 )$, for all $x,y\in\R^n$. $\mc{A}$ is  $\ell$-Lipschitz continuous in $\H_P$ if $\|\mc{A}(x)-\mc{A}(y)\|_P \leq \ell \|x-y\|_P$, for all $x,y\in\R^n$; if $\ell\leq 1$ ($\ell <1)$, $\mc{A}$ is nonexpansive (contractive) in $\H_P$. We omit the indication ``in $\H_P$" if $P=I$.
 $\proj _S^P:\R^n\rightarrow S $ is the Euclidean $P$-weighted projection onto a closed convex set $S\subseteq\R^n$, i.e. 
$\proj _S^P (x):=\argmin_{y\in S}\|y-x\|_P$; we omit the superscript if $P=I$.

\section{Mathematical setup}\label{sec:mathbackground}
\subsection{The game}
 We consider a set of  agents, $ \mc I:=\{ 1,\ldots,N \}$, where each agent $i\in \mc{I}$ shall choose its decision variable (i.e., strategy) $x_i$ from its local decision set $\textstyle \Omega_i \subseteq \R^{
n_i}$.  Let $x := \col( (x_i)_{i \in \mc I})  \in \Omega $ denote the stacked vector of all the agents' decisions, with $\textstyle \Omega := \Omega_1\times\dots\times\Omega_N\subseteq \R^n$ the overall action space and $n:=\textstyle \sum_{i\in \mc{I}} n_i$.
The goal of  agent $i \in \mc I$ is to minimize its objective function $J_i(x_i,x_{-i})$, which depends both on the local variable $x_i$ and on the decision variables of the other agents $x_{-i}:= \col( (x_j)_{j\in \mc I\backslash \{ i \} })$.
%
The game is then represented by the inter-dependent optimization problems
\begin{align} \label{eq:game}
\forall i \in \mc{I}:
\quad \underset{y_i \in \Omega_i}{\argmin}   \; J_i(y_i,x_{-i}).
\end{align} 
The technical problem we consider here is the distributed computation of a \gls{NE}, as formalized next.

\begin{definition}
	A collective strategy $x^{*}=\operatorname{col}\left((x_{i}^{*}\right)_{i \in \mathcal{I}})$ is a Nash equilibrium if, for all $i \in \mc{I}$,
	\[
	J_{i}\left(x_i^*, x_{-i}^{*}\right)\leq \inf \{J_{i}\left(y_{i}, x_{-i}^{*}\right) \mid (y_i,x_{-i}^*)\in\Omega \}. \QEDopenhereeqn 
	\]
\end{definition}

\smallskip
Next, we postulate  common regularity
assumptions for the constraint sets and cost functions  \cite[Ass.~1]{Pavel2018}, \cite[Ass.~1]{TatarenkoShiNedic_CDC2018}.

\begin{standing}\label{Ass:Convexity}
	For each $i\in \mathcal{I}$, the set $\Omega_i$ is non-empty, closed and convex;  $J_{i}$ is continuous and  $J_{i}\left(\cdot, x_{-i}\right)$ is convex and continuously differentiable for every $x_{-i}$.
	{\hfill $\square$} \end{standing}

Under Standing Assumption~\ref{Ass:Convexity}, a collective strategy $x^*$ is a \gls{NE} of the game in \eqref{eq:game} if and only if it is a solution of the variational inequality VI$(F,\Omega)$\footnote{Given a set $S\subseteq \R^m$ and a mapping $\psi:S\rightarrow \R^m$, the  VI$(\psi,S)$ is the problem of finding $\omega^*\in S$ such that $\langle \psi(\omega^*),\omega-\omega^*\rangle \geq 0$, for all $\omega \in S$.}  \cite[Prop.~1.4.2]{FacchineiPang2007},
where $F$ is the \emph{pseudo-gradient} mapping of the game:
\begin{align}
\label{eq:pseudo-gradient}
F(x):=\operatorname{col}\left( (\nabla _{\! \! x_i} J_i(x_i,x_{-i}))_{i\in\mathcal{I}}\right).
\end{align}
Equivalently, $x^*$ is a \gls{NE} if and only if  
\begin{align} \label{eq:NEinclusion}
\forall i \in \mc{I}: \quad x_i^*=\proj_{\Omega_i}(x_i^*-\beta_i \nabla_{\! \! x_i} J_i(x_i^*,x_{-i}^*)),
\end{align}
for arbitrary positive scalars $\beta_i$'s \cite[Prop.~1.5.8]{FacchineiPang2007}. A sufficient condition for the existence and uniqueness of a \gls{NE}  is the strong monotonicity of the pseudo-gradient \cite[Th. 2.3.3]{FacchineiPang2007}, as postulated next. This assumption has always been used for \gls{NE} seeking under partial-decision information with fixed step sizes, e.g.,   \cite[Ass.~2]{Pavel2018}, \cite[Ass.~4]{DePersisGrammatico2018}, \cite[Ass.~2]{TatarenkoShiNedic_CDC2018}.

\sloppy
\begin{standing}\label{Ass:StrMon}
	The pseudo-gradient mapping $F$ in \eqref{eq:pseudo-gradient}  is $\mu$-strongly monotone and $\ell_0$-Lipschitz continuous, for some $\mu$, $\ell_{0}>0$.
	\hfill $\square$
\end{standing}

\subsection{Network communication}
 The agents can exchange information with some neighbors over a directed communication network $\mathcal G(\mc{I},\mc{E})$.  The ordered pair $(i,j) $ belongs to the set of edges, $\mc{E}$, if and only if agent $i$ can receive information from agent $j$. 
We denote $W\in \R^{N\times N}$ the weighted adjacency matrix of $\mc{G}$ and $w_{i,j}:=[W]_{i,j}$, with $w_{i,j}>0$ if $(i,j)\in \mc{E}$, $w_{i,j}=0$ otherwise; 
$d_i=\textstyle \sum _{j=1}^N w_{i,j}$ and $\mc{N}_i=\{j\mid (i,j)\in \mc{E}_k\}$ the in-degree and  the set of in-neighbors of agent $i$, respectively. 
 

\begin{standing}\label{Ass:stronglyconnectedgraph}
	The communication graph $\mc{G}$ is strongly connected.  \hfill $\square$
\end{standing}

\begin{standing}\label{Ass:rowstochasticity}
	The adjacency matrix $W$ satisfies the following conditions:
	\begin{itemize}
		\item[(i)] \emph{Self-loops:} $w_{i,i}>0$ for all $i\in \mc{I}$;
		\item[(ii)] \emph{Row stochasticity:} $W_k\1_N=\1_N$.
		{\hfill $\square$}
	\end{itemize} 
\end{standing}

\begin{remark}
  Standing Assumption~\ref{Ass:rowstochasticity} can be fulfilled on any digraph, if the agents can access their own in-degree, by \emph{locally} assigning weights to the received information.  \hfill $\square$
\end{remark}

Under Standing Assumptions~\ref{Ass:stronglyconnectedgraph}-\ref{Ass:rowstochasticity}, by the \gls{PF} theorem, $W$ has a simple eigenvalue in $1$; all the other (complex) eigenvalues of $W$ have absolute value strictly smaller than $1$.  Besides, there exist a vector  $q=\col((q_i)_{i\in\mc{I}})$ such that
\begin{equation}
q\in \R^N_{>0}, \quad q^\top W=q^\top, \quad \1_N^\top q=1.
\end{equation}  
We call $q$ the (left) \emph{Perron-Frobenius eigenvector} of $W$. 
Let \begin{equation}
Q:=\diag((q_i)_{i\in\mc{I}}).
\end{equation}
Clearly, $Q\succ 0$. Unless $W$ is doubly stochastic, $W$ is not nonexpansive in $\H_I$, i.e., $\sigmamax{W}>1$. This is one of the main technical challenges to face when studying fixed-point iterations over directed graphs  \cite{SalehisadaghianiPavel2017_nondoublystochastic}.
To deal with this complication, it was shown in  \cite[Lemma~1]{Cenedese2019AsynchronousAT} that $W$ is nonexpansive (averaged, indeed) in $\H_Q$. Next, we provide an additional contractivity result, which we exploit later on.

\begin{lemma}\label{lem:Qcontractivity}
For any $y\in\R^N$, $\|W(y-\1_{N} q^\top y)\|_Q \leq \sigmabar \| y-\1_{N} q^\top y\|_Q$, where $\sigmabar:=\sigma_{N-1}(Q^{\scriptstyle \frac{1}{2}}WQ^{-\scriptstyle \frac{1}{2}})<1$. \hfill $\square$ 
\end{lemma}

If $W$ is also column  stochastic, Lemma~\ref{lem:Qcontractivity} holds with  $q=\textstyle \frac{1}{N}\1_N$ and  $Q=I_N$, and we recover a well-known property of doubly stochastic matrices \cite[Eq.~4]{Bianchi_LCSS2020}.

\begin{remark}
	\cite[Lemma~1]{XiMai_DistributedOptimizationonRowStochastic} states that there exist a norm and $\sigmabar>0$ such that the property in Lemma~\ref{lem:Qcontractivity} holds; instead, we  explicitly characterized  both the norm and $\sigmabar$, which proves very advantageous in our analysis, see §\ref{sec:distributedGNE}. \hfill $\square$
\end{remark}

\subsection{Partial-decision information scenario}
 In our setup,  agent $i\in\mc{I}$ can only access its own feasible set $\Omega_i$  and an analytic expression of its own cost function $J_i$. However, the agents cannot evaluate the actual value of the cost $J_i(x_i,x_{-i})$ (or the partial derivative $\grad_{\!\!x_i}J_i(x_i,x_{-i})$), since they cannot access the strategies of all the competitors $x_{-i}$. 
Instead, the  agents only rely on the information exchanged locally with their neighbors over the communication graph $\mc{G}$. To cope with the lack of knowledge, the general assumption for this partial-decision information scenario is that each agent keeps an estimate of all other agents' actions \cite{Pavel2018}, \cite{Koshal_Nedic_Shanbag_2016}, \cite{DePersisGrammatico2018}. Then, the agents aim at reconstructing the actual values, based on the data received from their neighbors. 
We denote $\x_{i}=\operatorname{col}((\x_{i,j})_{j\in \mc{I}})\in \R^{n}$,  where $\x_{i,i}:=x_i$ and $\x_{i,j}$ is agent $i$'s estimate of agent $j$'s action, for all $j\neq i$; $\x_{j,-i}=\col((\x_{j,l})_{l\in\mc{I}\backslash \{ i \}})$;  $\x=\col((\x_i)_{i\in\mc{I}})$. As in \cite[Eq.13-14]{Pavel2018}, we define
	\begin{align}
\mathcal{R}_{i}:=&\left[ \begin{array}{lll}{{0}_{n_{i} \times n_{<i}}} & {I_{n_{i}}} & {\0_{n_{i} \times n_{>i}}}\end{array}\right], 
\end{align}
where $n_{<i}:=\sum_{j<i,j \in \mathcal{I}} n_{j}$, $n_{>i}:=\sum_{j>i, j \in \mathcal{I}} n_{j}$. 
In simple terms,  
$\mathcal R _i$ selects the $i$-th $n_i$-dimensional component from an $n$-dimensional vector, i.e.,  $\mathcal{R}_{i} \x_{i}=\x_{i,i}=x_i$.
We denote by $\mathcal{R}:=\operatorname{diag}\left((\mathcal{R}_{i})_{i \in \mathcal{I}}\right)$;
 thus, we have $x=\mathcal{R} \x$. Moreover, we define  the  \textit{extended pseudo-gradient} mapping  $\bs{F}$ as
\begin{align}
\label{eq:extended_pseudo-gradient}
\bs{F}(\x):=\operatorname{col}\left((\nabla_{\! \! x_{i}} J_{i}\left(x_{i}, \x_{i,-i}\right))_{i \in \mathcal{I}}\right).
\end{align}
\begin{lemma}[\!\!{\cite[Lemma 3]{Bianchi_ECC20_ctGNE}}]\label{lem:LipschitzExtPseudo}
	The  mapping $\bs{F}$ in \eqref{eq:extended_pseudo-gradient} is $\ell$-Lipschitz continuous, for some  $\ell\in [\mu, \ell_0]$: for any $\x,\y\in\R^{Nn}$, $\|\F(\x)-\F(\y)\|\leq \ell\|\x-\y\|$.
	{\hfill $\square$} \end{lemma}

We remark that in \eqref{eq:extended_pseudo-gradient}, each agent $i$ evaluates its partial gradients  $\grad_{\!\!x_i}J_{i}\left(x_{i}, \x_{i,-i}\right)$  on the local estimate $\x_{i,-i}$, not on the actual strategies $x_{-i}$. Only when the estimates of all the agents coincide with the actual value, i.e., $\x =\1_N\otimes x$, we have that $\F(\x)=F(x)$. As a consequence, the mapping $\mc{R}^\top\F$ is not monotone, not even under strong monotonicity of the game mapping  $F$ in Standing Assumption~\ref{Ass:StrMon}. Indeed, the loss of monotonicity is the main technical difficulty arising in the partial-decision information scenario \cite{Pavel2018}, \cite{TatarenkoShiNedic_CDC2018}. 

\section{Fully-distributed  Nash equilibrium seeking}\label{sec:distributedGNE} In this section, we present a pseudo-gradient method (along with some variants) to seek a \gls{NE} in a  fully-distributed way. Before going into details, we need some  definitions. Let
\begin{equation}\label{eq:Qbs}
\Qbar:=\diag((q_i I_{n_i})_{i\in\mc{I}}), \quad  \Q:=Q\otimes I_n.
\end{equation}
We define the consensus subspace as $\bs{E}=\{\y \in\R^{Nn} | \y=\1_N\otimes y, y\in\R^n \}$ and its orthogonal complement in $\H_{\Q}$ as $\bs{E}^Q_\perp=\{ \y\in\R^{Nn} | (q \otimes I_n)^\top \y=\0_n\}$. Thus, any vector of estimates $\x\in\R^{Nn}$ can be written as $\x=\x_{\scriptscriptstyle \parallel}+\x_{\! \scriptscriptstyle \perp}$, where $\x_{\scriptscriptstyle \parallel}=\proj^{\Q}_{\bs E}(\x)=(\1_N q^\top\otimes I_n) \x$, $\x_{\! \scriptscriptstyle \perp}=\proj^{\Q}_{\bs E^Q_{\perp}}(\x)$, and it holds that $\langle \x_{\scriptscriptstyle \parallel},\x_{\! \scriptscriptstyle \perp}\rangle_{\Q}=0$. Clearly, if the estimates of the agents $\x \in \bs{E}$, then $\x_i=x$ for all $i\in \mc{I}$, namely the estimate of each agent coincides with the actual collective strategy $x$.

\subsection{Case 1: Known $q$ and $\bar{\sigma}$}\label{sec:math:A}

Our basic fully-distributed \gls{NE} seeking algorithm is summarized in Algorithm~\ref{algo:1}, where $\alpha$ is a fixed step size.  Each agent update its estimates according to consensus
dynamics, then its strategy via a projected pseudo-gradient step. We remark
that each agent computes the partial gradient of its cost in
its local estimate, not on the actual joint strategy $x$. 

Compared to similar pseudo-gradient dynamics proposed in the literature \cite{TatarenkoShiNedic_CDC2018}, \cite{Bianchi_LCSS2020}, the novelty of  Algorithm~\ref{algo:1} is that the cost related components $\nabla_{\! \! x_i} J_i$ are weighted  by the reciprocal of the elements $q_i$ of the \gls{PF} eigenvector. This operation enables convergence on row stochastic graphs, and in fact it is not necessary for doubly stochastic graphs, for which $q=\1$.
The idea behind this key modification is that $(W-\1 q^\top)$ is contractive in $\H_Q$, while the game-mapping $F$ is strongly monotone in $\H_I$; instead, we would like both properties to hold in the same space. Division by the \gls{PF} eigenvector achieves this goal, as we show next. Let 
\begin{equation}
\fq(x)=:{\Qbar}^{-1}F(x),  \quad \FQ(\x):={\Qbar}^{-1}\F(\x).
\end{equation}

\begin{lemma}\label{lem:Qlipmon}
	$\fq$ is $\barmu$-strongly monotone  and $\barell_0$-Lipschitz continous in $\H_{\Qbar}$,  for some $\barmu, \barell_0 >0$; 
	$\FQ$ is $\barell$-Lipschitz continuous from  $\H_{\Q}$ to $\H_{\Qbar}$, for some $\barell>0$, i.e., for any $\x,\y \in \R^{Nn}$, $\|\FQ(\x)-\FQ(\y)\|_{\Qbar}\ \leq \barell\|\x-\y\|_{\Q}$. \hfill $\square$
\end{lemma}


\begin{remark}
	Lemmas~\ref{lem:Qcontractivity} and \ref{lem:Qlipmon} provide a general, operator-theoretic interpretation of the approach in \cite{XiMai_DistributedOptimizationonRowStochastic}, where a similar  technique is used in the context of  distributed optimization. \hfill $\square$
\end{remark}
 
 \newlength{\textfloatsepsave} \setlength{\textfloatsepsave}{\textfloatsep} \setlength{\textfloatsep}{10pt} 
\begin{algorithm}[t]\caption{} \label{algo:1}
	\textbf{Initialization}: $\forall i\in \mc{I}$, set $x_i^0\in \Omega_i$, $\bs{x}_{i,-i}^0\in \R^{N-n_i}$. 
	\\ 
	\textbf{Iterate until convergence:} each agent $i\in\mc{I}$ does:
	\begin{align*}
	\hat{\x}_i^k &= \textstyle \sum_{j\in\mc{N}_i} w_{i,j}^k  \x_{j}^k \\
		x_i^{k+1} &=\proj_{\Omega_i} (\xhb_{i,i}^k-\textstyle  \frac{\alpha}{q_i} \nabla_{ \! \! x_i} J_i (\hat{\x}^k_i))\\
		\x_{i,-i}^{k+1} &=\hat{\x}^k_{i,-i}.
		\end{align*}
		\vspace{-1.3em}
\end{algorithm}

In compact form, Algorithm~\ref{algo:1} reads as 
\begin{equation}\label{eq:algo1:compact}
\x^{k+1}=\proj_{\bs{\Omega}}(\mc{F}(\x^k)),
\end{equation}
where $\boldsymbol{\Omega}:=\{\boldsymbol{x} \in \mathbb{R}^{N n} \mid \mathcal{R} \boldsymbol{x} \in \Omega \}$, $\bs{W}:=W \otimes I_n$ and
\begin{equation}\label{eq:opF}
\mc{F}(\x):=\W \x-\alpha\mc{R}^\top \FQ(\W\x).
\end{equation}
 The following Lemma shows a contractivity property of the operator $\mc{F}$ and represents the cornerstone we use to prove convergence of our \gls{NE}
 seeking schemes. The result is based on the strong monotonicity of $\fq$ in $\H_{\bar{Q}}$ and on Lemma~\ref{lem:Qcontractivity}. 

\begin{lemma}\label{lem:contractivityF} 
	Let
	\begin{align}\label{eq:M}
	M_{\alpha}  \! :=  \!
	\begin{bmatrix}
	1 \! - \! 2\alpha \barmu\eigmin{Q} \! +\! \alpha^2 \barell^2 &  (2\alpha\barell)\sigmabar\\
	 (2\alpha\barell)\sigmabar & ( 1\! + \! 2\alpha\barell \! +\! \alpha^2\barell^2 ) {\bar{\sigma}}^2 
	\end{bmatrix} 
    \end{align}
	If the step size $\alpha>0$ is chosen such that 
	\begin{equation}\label{eq:C1}
	\rho_\alpha:=\uplambda_{\max}(M_{\alpha})=\|M_\alpha\|<1,
	\end{equation}
	then the operator $\mc{F}$ in \eqref{eq:opF} is $\sqrt{\rho_\alpha}$-restricted contractive in $\H_{\Q}$ with respect to the consensus subspace $\bs{E}$, i.e., for any $\x\in\R^{Nn}$, $\y \in\bs{E}$, it holds that 
	$\| \mc{F}(\x)-\mc{F}(\y) \|_{\Q} \leq  \sqrt{\rho_\alpha} \| \x-\y\|_{\Q}. $ \hfill \QEDopen
\end{lemma}

\begin{remark}
	The condition in \eqref{eq:C1} can always be satisfied by choosing $\alpha$  small enough; an  explicit upper bound can be obtained as in \cite[Lemma~2]{Bianchi_LCSS2020}.  \hfill $\square$
\end{remark}

 
\begin{theorem}\label{th:main1}
		Let $\alpha>0$ satisfy the condition in \eqref{eq:C1}. Then, for any initial condition, the sequence $(\x^k)_{k\in\N}$ generated by  Algorithm~\ref{algo:1} converges to  $\x^*=\1_N\otimes x^*$, where $x^*$ is the  \gls{NE} of the game in \eqref{eq:game}, with linear rate: for all $\k$, 
	\[
	\| \x^k -\x^*\|_{\Q} \leq \left(\sqrt{\rho_\alpha}\right)^{\; k} \| \x^0-\x^*\|_{\Q}.
	\QEDopenhereeqn
	\]
\end{theorem}
\smallskip

\begin{proof}
	By \eqref{eq:NEinclusion}, we  infer that $x^*$ is the \gls{NE} if and only if $x^*=\proj_\Omega(x^*-\alpha Q^{-1} F(x^*))$. Together with $\W\x^*=\x^*$ and $\F(\x^*)=F(x^*)$, this implies that $\x^*$ is a fixed point for the iteration in \eqref{eq:algo1:compact}. Therefore we can write
	\[\begin{aligned}
     \| \x^{k+1}-\x^*\|_{\Q} &=\| \proj_{\bs{\Omega}} (\mc{F} (\x^{k}))-\proj_{\bs{\Omega}}(\mc{F}(\x^*)) \|_{\Q}\\ 
	& =\| \proj^{\Q}_{\bs{\Omega}} (\mc{F} (\x^{k}))-\proj^{\Q}_{\bs{\Omega}}(\mc{F}(\x^*)) \|_{\Q} \\
	& \leq \| \mc{F} (\x^{k })-\mc{F}(\x^*)\|_{\Q}\leq \sqrt{\rho_\alpha} \|\x^k-\x^*\|_{\Q},
	\end{aligned}
	\]
	where the second equality follows by $\Q=Q\otimes I_n$ and the definition of $\bs{\Omega}$ (note that $\proj _{\Omega_i}^{q_i I_{n_i}}=\proj_{\Omega_i}$), the first inequality follows by nonexpansiveness of the projection \cite[Prop. 4.16]{Bauschke2017}, and the second inequality by Lemma~\ref{lem:contractivityF}. 
\end{proof}

We note that Algorithm~\ref{algo:1} requires a priori knowledge of the   communication graph $\mc{G}$, both to compute the \gls{PF} eigenvector $q$ and to tune the step size $\alpha$. In the next subsection, we relax this hypothesis.

\subsection{Case 2: Online computation of $q$}\label{sec:math:B}
When the \gls{PF} eigenvalue $q$ is not known in advance, it can be computed online in a distributed fashion. The procedure is illustrated in Algorithm~\ref{algo:2}. Each agent $i\in\mc{I}$ keeps an extra variable $\hat{q}_i=\col((\hat{q}_{i,j})_{j\in\mc{I}})$, which is an estimate of $q$, initialized as the $i$-th vector of the canonical basis $e_i^N\in\R^N$. 

Notably, each estimate $\hat{q}_i$ converges to the real value $q$. In fact, the updates in Algorithm~\ref{algo:2} can be written compactly as 
\begin{equation}\label{eq:PGFestimateupdate}
 \hat{\bs{q}}^{k+1}=(W\otimes I_N)\hat{\bs{q}}^k,
\end{equation}
where $\hat{\bs{q}}:=\col((\hat{q}_i)_{i\in\mc{I}})$. Therefore, by the \gls{PF} theorem (and by Standing Assumptions~\ref{Ass:stronglyconnectedgraph}-\ref{Ass:rowstochasticity}), $\hat{\bs{q}}^k$ converges linearly to $(\1_N q^\top \otimes I_N) \hat{\bs{q}}^0=\1_N \otimes q$.
In particular, $\hat{q}^k_{i,i} \rightarrow q_i$. Also, $\hat{q}^k_{i,i}>0$ for all $k\geq 0$, since $\hat{q}^0_{i,i}>0$ and $W$ is nonnegative with positive diagonal. As such, Algorithm~\ref{algo:2} is always well defined. We first show its convergence for a fixed step size.

\begin{algorithm}[t]\caption{} \label{algo:2}
	\textbf{Initialization}: $\forall i\in \mc{I}$, set $x_i^0 \hspace{-0.1em}\in \hspace{-0.1em} \Omega_i$, $\bs{x}_{i,-i}^0\in \R^{N-n_i}$, $\hat{q}_i^0=e_i^N\!$.
	\\ 
	\textbf{Iterate until convergence:} each agent $i\in\mc{I}$ does:
		\begin{align*}
		\hat{q}_i^{k+1} &= \textstyle  \sum_{j\in\mc{N}_i} w_{i,j}^k \hat{q}_j^k \\
		x_i^{k+1}&=\proj_{\Omega_i} (\xhb_{i,i}^k- \alpha^k(\hat{q}_{i,i}^k)^{-1} \nabla_{\!\! x_i} J_i (\hat{\x}^k_i))\\
		\x_{i,-i}^{k+1}&=\hat{\x}^k_{i,-i} \qquad
		\hat{\x}_i^k  = \textstyle \sum_{j\in\mc{N}_i} w_{i,j}^k  \x_{j}^k.
		\end{align*}
		\vspace{-1.3em}
\end{algorithm}

\begin{theorem}\label{th:main2}
	Let $\alpha>0$ satisfy  the condition in \eqref{eq:C1}, and $\alpha^k=\alpha$   $\forall \k$.  Then, for any initial condition, the sequence $(\x^k)_{\k}$ generated by Algorithm~\ref{algo:2} converges to  $\x^* =\1_N\otimes x^*$, where $x^*$ is the  \gls{NE} of the game in \eqref{eq:game}, with linear rate:  for any $\epsilon>0$,  there exists $K>0$ such that, for all $\k$,
	\[
	\| \x^k -\x^*\|_{\Q} \leq K \left(\sqrt{\rho_\alpha}+\epsilon\right)^{\; k} \| \x^0-\x^*\|_{\Q}. \QEDopenhereeqn
	\]
\end{theorem}
\smallskip

While in Algorithm~\ref{algo:2} the \gls{PF} eigenvector is estimated online,  the upper bound on  $\alpha$ in Theorem~\ref{th:main2} is still a function of the network parameter  $\sigmabar$, which can be difficult to compute distributedly. Upper/lower bounds might be available for some classes of networks,
e.g., unweighted graphs. This is analogous to \cite[Th.~2]{XiMai_DistributedOptimizationonRowStochastic}, where $q$ is computed online, but the step size depends on global, not easily accessible, information. In fact, this notion of fixed but small-enough step sizes is not uncommon in distributed algorithms literature.
%
%

 When estimating a step $\alpha$ that satisfies \eqref{eq:C1} is impossible, convergence to a \gls{NE} can still be guaranteed by allowing for diminishing step sizes.
 In this case, also the information  on the game (i.e., Lipschitz and monotonicity constants of the pseudo-gradient) is not needed for the tuning.

\begin{theorem}\label{th:main3}
	Let $(\alpha^k)_{\k}$ be a positive nonincreasing sequence such that  $\textstyle \sum_{k\in \N} \alpha^k=\infty$ and $\lim_{k\to\infty} \alpha^k=0$. Then, for any initial condition, the sequence $(\x^k)_{k\in\N}$ generated by Algorithm~\ref{algo:2} converges to  $\x^*=\1_N\otimes x^*$, where $x^*$ is the \gls{NE} of the game in \eqref{eq:game}. \hfill $\square$ 
\end{theorem}

\section{Numerical example: A Nash-Cournot game}\label{sec:numerics}

\begin{figure}[t]
 	\centering
 	\includegraphics[width=1\columnwidth]{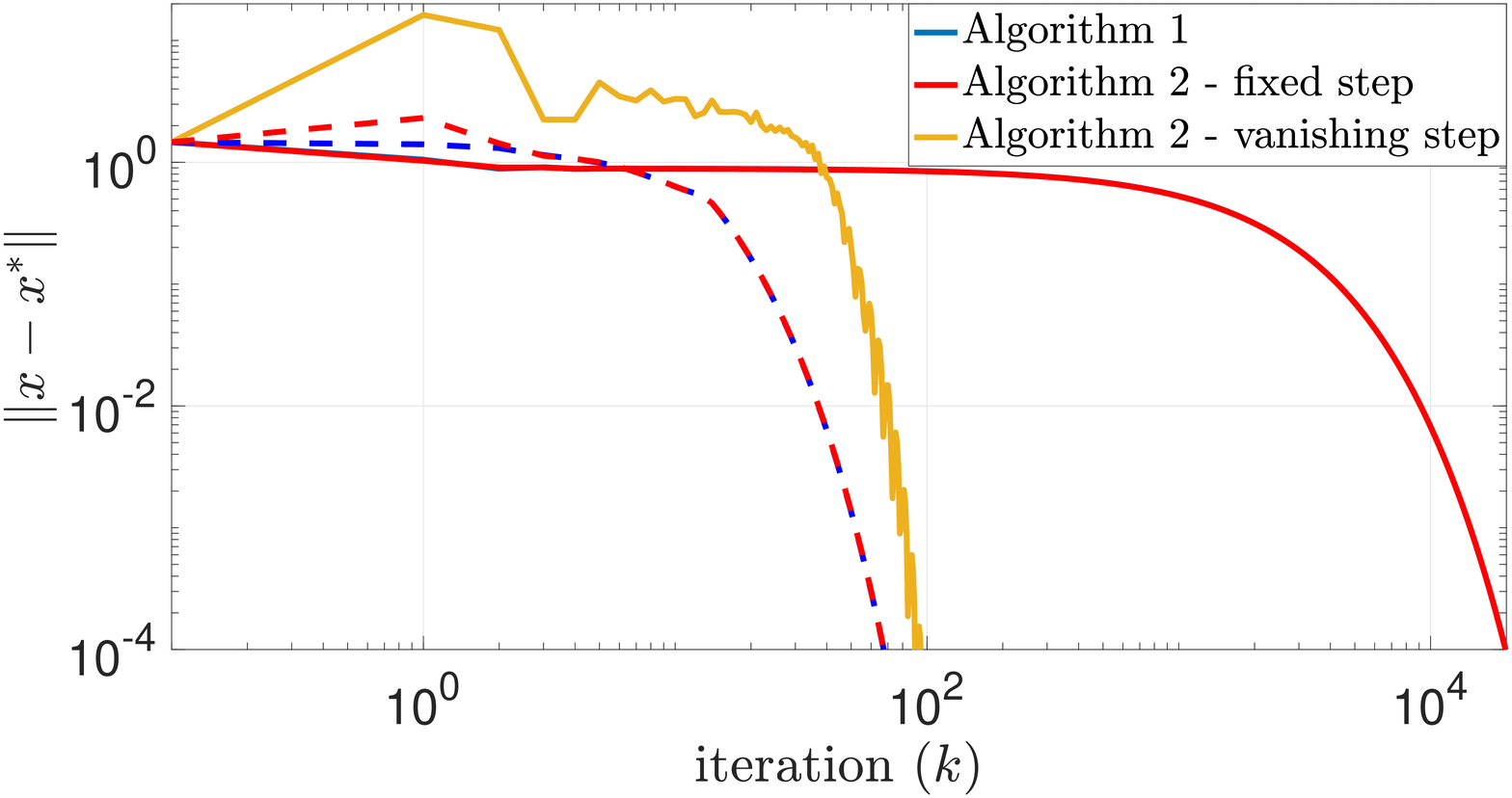}
 	\caption{Distance from the Nash equilibrium, with the step sizes that ensure convergence (solid lines) and with a fixed step size chosen $400$ times bigger than the theoretical upper bound (dashed lines).}\label{fig:1}
 \end{figure}

We consider the  Cournot competition model in \cite[§6]{Pavel2018}. $N$ firms produce an uniform commodity that is sold to $m$ markets. Each firm $i\in\mc{I}=\{1,\dots,N\}$ is allowed to participate in $n_i\leq m$ of the markets; its decision variable is the vector $x_i\in\R^{n_i}$ of quantities of product to be delivered to each of the $n_i$ markets, bounded by the local constraints $\0_{n_i}\leq x_i\leq X_i$. Let $A_i\in\R^{m\times n_i}$  such that $[A_i]_{k,j}=1$ if $[x_i]_j$ is the amount of commodity sent to the $k$-th market by agent $i$, $[A_i]_{k,j}=0$ otherwise, for all $j=1,\dots, n_i$, $k=1,...,m$.
Hence,   $Ax=\textstyle \sum_{i=1}^N A_ix_i \in \R^m$, where $A:=[A_1 \dots A_N]$, are the quantities of product delivered to each market.  Firm $i$ aims at maximizing its profit, i.e.,  minimizing the cost
$J_i(x_i,x_{-i})=c_i(x_i)-p(Ax)^\top A_ix_i.$
Here, $c_i(x_i)=x_i ^\top Q_i x_i+q_i^\top x_i$ is firm $i$'s production cost, with $Q_i\succ 0$; $p:\R^m\rightarrow \R^m$ associates to each market a price that depends on the amount of product delivered to that market. Specifically, for $k=1,\dots,m$, $[p(x)]_k=\bar P_k$ -$\chi_k [Ax]_k$, where $\bar P_k, \chi_k>0$. 
We set $N=20$, $m=7$. The market structure (i.e., which firms participate in each market) is defined as in \cite[Fig.~1]{Pavel2018}. Therefore, $x=\col((x_i))_{i\in\mc{I}})\in\R^n$ and $n=32$.  We select randomly with uniform distribution $r_k$ in $[1,2]$, $Q_i$ diagonal with diagonal entries in $[14,16]$, $q_i$ with elements in $[1,2]$, $\bar{P}_k$ in $[10,20]$, $\chi_k$ in $[1,3]$, $X_i$  in $[5,10]$, for all $i\in \mc{I}$, $k=1,\dots,m$.
This setup satisfies Standing Assumptions \ref{Ass:Convexity}-\ref{Ass:StrMon}  \cite[§6]{Pavel2018}. 
The firms communicate  over a randomly generated strongly connected row stochastic directed network, but cannot access the production of all the competitors. 
We set $\alpha\approx 3\times 10^{-5}$ to satisfy the condition in \eqref{eq:C1}. We compare the performance of Algorithms~\ref{algo:1} and Algorithm \ref{algo:2}, the latter both with a fixed  ($\alpha^k=\alpha$) and  vanishing step size ($\alpha^k=\frac{1}{k+1}$), in figure~\ref{fig:1}. Due to the small $\alpha$,  the schemes with fixed step are almost indistinguishable, and diminishing step sizes result in faster convergence. The good performance obtained with vanishing step suggests that the choice of $\alpha$ is quite conservative. Indeed, Algorithms~\ref{algo:1}-\ref{algo:2} still converge, and much faster,  with a fixed step size  $400$ times larger than its theoretical upper bound (dashed lines).

\section{Conclusion}\label{sec:conclusion}
Certain properties of doubly stochastic matrices  carry on to row stochastic matrices, but in a different Hilbert space, weighted by their left Perron-Frobenius eigenvector. We exploited one such contractivity property to solve, in a fully-distributed way, Nash equilibrium problems over directed networks. Any requirement for global knowledge of the graph and of the game mapping can be avoided in the case of vanishing step sizes. 

The extension of our results to generalized games, where the agents share some common constraints, is left as future research.
It would be also  valuable to relax our connectivity and monotonicity assumptions, namely allowing for  jointly  connected networks and (strictly) monotone game mappings.

\appendix
\subsubsection{Proof of Lemma~\ref{lem:Qcontractivity}}\label{app:lem:Qcontractivity}
\noindent By $W\1_N=\1_N$, it suffices to show that $\|W-\1_N q^\top\|_Q=\sigmabar<1$. Let $p:=\col((\sqrt{ q_i})_{i\in\mc{I}})$. Then,
\begin{align*}
& \ \quad \|W-\1_N  q^\top\|^2_Q =\|Q^{\scriptstyle \frac{1}{2}}(W-\1_N  q^\top)Q^{-\scriptstyle \frac{1}{2}}\|^2 
\\
& = \eigmax{ (Q^{\scriptstyle \frac{1}{2}}WQ^{-\scriptstyle \frac{1}{2}}-pp^\top)^\top(Q^{\scriptstyle \frac{1}{2}}WQ^{-\scriptstyle \frac{1}{2}}-pp^\top)} 
\\ & \overset{(\textnormal{a})}{=} \eigmax{Q^{-\scriptstyle \frac{1}{2}}W^\top Q W Q^{-\scriptstyle \frac{1}{2}}-pp^\top}:=\eigmax{M-pp^\top},
\end{align*}
where in (a) we used $p^\top p=1$, and $M=Q^{-\scriptstyle \frac{1}{2}}W^\top Q W Q^{-\scriptstyle \frac{1}{2}}$.
Since $M$ is symmetric and $ M p=p$, $M$ has a basis of eigenvectors, say 
$\{v_1,\dots, v_{N-1},p \}$, with associate eigenvalues $\{s_1,\dots,s_{N-1},1\}$. By orthogonality and $p^\top p=1$, it follows that the eigenvalues of $M-pp^\top$ are $\{s_1,\dots,s_{N-1},0\}$, with associate eigenvectors $\{v_1,\dots,v_{N-1},p\}$. Since $M\succeq 0$, it suffices to show that $s_i<1$, for $i=1, \dots, N-1$. Seeking a contradiction, let ${j}\in \{1, \dots, N-1 \}$ such that $s_{j}\geq 1$, and  $\bar{v}:=Q^{\scriptstyle -\frac{1}{2}}v_j$. Thus, we have $\|W\bar{v}\|^2_Q=v_j^\top Q^{\scriptstyle -\frac{1}{2}} W^\top Q W Q^{\scriptstyle -\frac{1}{2}}  v_j=v_j^\top M v_j \geq v_j^\top v_j = \bar{v}Q\bar v=\|\bar{v}\|_Q^2$. 
By \cite[Lemma~1]{Cenedese2019AsynchronousAT}, it also holds, for some $\gamma>0$, for any $y\in \R^N$, that $
\|Wy\|_Q\leq\|y\|_Q-\gamma \|(I_N-W)y\|_Q
$. Hence, by Standing Assumption~\ref{Ass:stronglyconnectedgraph}, it must hold that $\bar{v}=\beta 1_N$, for some $ \beta\neq 0$. Equivalently, $v_j=\beta p$. This is a contradiction, since $p$ and $v_j$ must be  orthogonal. The conclusion follows with $\bar{\sigma}=\sqrt{\eig_{N-1}(M)}$. \hfill $\blacksquare$

\smallskip
\subsubsection{Proof of Lemma~\ref{lem:Qlipmon}}\label{app:lem:Qlipmon}
\noindent For any $x,y\in\R^n$ it holds that 
$ \langle {\Qbar}^{-1}(F(x)-F(y)), x-y\rangle_{{\Qbar}} =  \langle F(x)-F(y), x-y\rangle  
\geq \mu \|x-y\|^2 \geq \textstyle \frac{\mu}{\eigmax{{{\Qbar}}}}\|x-y\|^2_{{{\Qbar}}}$, and that $\| {{\Qbar}}^{-1}(F(x)-F(y))\|^2_{\Qbar}=\|F(x)-F(y)\|^2_{{\Qbar}^{-1}}\leq \textstyle\frac{\eigmax{{\Qbar}^{-1}}}{\eigmin{{\Qbar}}}\ell_0^2\|x-y\|_{\Qbar}^2$.
Analogously, by Lemma~\ref{lem:LipschitzExtPseudo}, it holds that, for any $\x,\y\in\R^{Nn}$, $\| {{\Qbar}}^{-1}(\F(\x)-\F(\y))\|^2_{\Qbar}\leq \textstyle\frac{\eigmax{{\Qbar}^{-1}}}{\eigmin{\Q}}\ell^2\|\x-\y\|_{\Q}^2$. \hfill \QED

\smallskip

\subsubsection{Proof of Lemma~\ref{lem:contractivityF}}\label{app:lem:contractivityF}
We  use the shorthand notation $\FQ\x$ and $\fq x$ in place of $\FQ(\x)$ and $\fq(x)$. Let $\x\in\R^{Nn}$, $\y=\1\otimes y \in\bs{E}$, and  $\xhb:=\bs{W}\x=\xhb_{\scriptscriptstyle \parallel }+\xhb_{\! \scriptscriptstyle \perp}=\1_N\otimes \hat{x}_{\scriptscriptstyle \parallel}+\xhb_{\! \scriptscriptstyle \perp}\in\R^{Nn }$, with $\xhb_{\! \scriptscriptstyle \perp}\in E^Q_{\perp}$.
Thus, we have
%
%
\allowdisplaybreaks
\begin{align*} 
&\quad \  \| \mc{F}(\x)-\mc{F}(\y) \|^2_{\Q}
\nonumber 
& 
\\
\nonumber
& =  \|  (\xhb-\alpha\mc{R}^\top \FQ\xhb)-  (\y-\alpha\mc{R}^\top \FQ\y)\|_{\Q}^2 
%
\\ 
&\begin{aligned}[c]
&= \|\xhb_{\scriptscriptstyle \parallel}-\y\|_{\Q}^2+\|\xhb_{\! \scriptscriptstyle \perp} \|^2_{\Q} +\alpha^2  \|\mc{R}^\top(\FQ\xhb
-\FQ \y)\|_{\Q}^2 \\ 
& \qquad-2\alpha \langle \xhb_{\! \scriptscriptstyle \perp}, \mc{R}^\top(\FQ \xhb-\FQ \y)\rangle _{\Q} 
\\
& \qquad-2\alpha\langle \xhb_{\scriptscriptstyle \parallel}-\y, \mc{R}^\top (\FQ\xhb-\FQ \xhb_{\scriptscriptstyle \parallel})\rangle _{\Q} \\
& \qquad-2\alpha\langle \xhb_{\scriptscriptstyle \parallel}-\y,\mc{R}^\top(\FQ\xhb_{\scriptscriptstyle \parallel}-\FQ\y) \rangle _{\Q}
\end{aligned} 
\numberthis \label{eq:addends}
\\
\nonumber
& \leq \begin{aligned}[t]
& \|\xhb_{\scriptscriptstyle \parallel}-\y\|_{\Q}^2+\|\xhb_{\! \scriptscriptstyle \perp}\|_{\Q}^2+\alpha^2\barell^2(\|\xhb_{\! \scriptscriptstyle \perp}\|_{{\Q}}^2+ \|\xhb_{\scriptscriptstyle \parallel}-\y\|_{\Q}^2)
\\ &
+2\alpha\barell\|\xhb_{\! \scriptscriptstyle \perp}\| (\|\xhb_{\! \scriptscriptstyle \perp}\|_{\Q}+\|\xhb_{\scriptscriptstyle \parallel}-\y\|_{\Q})
\\& +\textstyle {2\alpha\barell }\|\xhb_{\scriptscriptstyle \parallel}-\y\|_{\Q}\|\xhb_{\! \scriptscriptstyle \perp}\|_{\Q}
\textstyle  -{2\alpha\barmu}\eigmin{Q}\|\xhb_{\scriptscriptstyle \parallel}-\y\|^2_{\Q},
\end{aligned}
\end{align*}
and to bound the addends in \eqref{eq:addends} we used:
\begin{itemize}[leftmargin=*]
	\item 3\textsuperscript{rd},  4\textsuperscript{th}, 5\textsuperscript{th} terms:   Lipschitz continuity of $\FQ$, the Cauchy-Schwartz inequality,  $\| \mc{R}^\top v \|_{\Q}=\|v\|_{\bar{Q}}$ for any $v\in\R^n$, 
	 $\|\xhb-\y\|_{\Q}^2=\|\xhb_{\scriptscriptstyle \parallel }-\y\|_{\Q}^2+\|\xhb_{\! \scriptscriptstyle \perp}\|_{\Q}^2$ (by orthogonality);
	\item 7\textsuperscript{th} term: $\langle \xhb_{\scriptscriptstyle \parallel}-\y,\mc{R}^\top(\FQ\xhb_{\scriptscriptstyle \parallel}-\FQ\y)\rangle _{\Q}=\langle \hat{x}_{\scriptscriptstyle \parallel}-y, \fq \hat{x}_{\scriptscriptstyle \parallel}-\fq  y\rangle_{{\Qbar}}\geq \barmu \|\hat{x}_{\scriptscriptstyle \parallel}-y\|^2_{{\Qbar}}\geq \mubar \eigmin{{\Qbar}}\|\hat{x}_{\scriptscriptstyle \parallel}-y\|^2=\mubar \eigmin{{\Qbar}}\|\xhb_{\scriptscriptstyle \parallel}-\y\|^2_{\Q}$, and the last equality follows since $\xhb_{\scriptscriptstyle \parallel},\y \in\bs{E}$, $\Q=Q\otimes I_n$ and  $\1_N^\top q=1$.
\end{itemize}
Besides, for every $\x=\x_{\scriptscriptstyle \parallel}+\x_{\! \scriptscriptstyle \perp}\in \R^{Nn}$, with $\x_{\scriptscriptstyle \parallel}\in\bs{E}$ and $\x_{\! \scriptscriptstyle \perp}\in \bs{E}^Q_{\! \perp}$, it holds that $ \xhb= \W\x=\x_{\scriptscriptstyle \parallel}+\W \x_{\! \scriptscriptstyle \perp}$, where $\W \x_{\! \scriptscriptstyle \perp} \in \bs{E}^Q_{\! \scriptscriptstyle \perp}$ (since  $(q\otimes I_n)^\top \W\x_{\! \scriptscriptstyle \perp}=(q\otimes I_n)^\top \x_{\! \scriptscriptstyle \perp}=\0_n$,  by definition of $\W$ and $q$). Consequently,  by Lemma~\ref{lem:Qcontractivity} and by  $\x_{\! \scriptscriptstyle \perp}=(I_{Nn}-\1_N q^\top\otimes I_n) \x$, we have 
$ \| \xhb_{\! \scriptscriptstyle \perp}\|_{\Q}=\| \W  \x_{\! \scriptscriptstyle \perp} \|_{\Q}\leq {\bar{\sigma}} \|\x_{\! \scriptscriptstyle \perp}\|_{\Q}.$
Therefore, we can finally write
\[
\begin{aligned}[b]
 \| \mc{F}(\x)-\mc{F}(\y)) \|^2_{\Q}
& \leq \left[\begin{matrix}
\|\x_{\scriptscriptstyle \parallel} -\y\|_{\Q} \\\nonumber \| \x_{\! \scriptscriptstyle \perp} \|_{\Q}
\end{matrix}\right]^\top 
M_\alpha 
\left[\begin{matrix}
\|\x_{\scriptscriptstyle \parallel} -\y\|_{\Q} \\\nonumber \|\x_{\! \scriptscriptstyle \perp}\|_{\Q}
\end{matrix}\right] 
\\
&  \leq \uplambda_{\max}(M_\alpha)(\|\texttt{}\x_{\scriptscriptstyle \parallel} -\y\|_{\Q}^2+\|\x_{\! \scriptscriptstyle \perp}\|_{\Q}^2) 
\\ 
&=\uplambda_{\max}(M_\alpha)\|\x-\y\|_{\Q}^2.   && \QED
\end{aligned}
\]

\subsubsection{Proof of Theorem~\ref{th:main2}}\label{app:th:main2}
We recast  Algorithm~\ref{algo:2} as 
\[
\x^{k+1}=\proj_{\bs{\Omega}} (\hat{\mc{F}}^k(\x^k)),
\]
where $\hat{\mc{F}}^k(\x^k) :=\W \x^k-\alpha\mc{R}^\top ({{\Qbar}}+\tilde{Q}^k)^{-1} \F(\W\x^k),$ and  $\tilde{Q}^k=\diag(((\hat{q}_{i,i}^k-q_i)I_{n_i})_{i\in\mc{I}})$. 
We  noted in §\ref{sec:math:B} that
 $(\bar{Q}+\tilde{Q}^k)=\diag((\hat{q}^k_{i,i}I_{n_i})_{i\in\mc{I}})\succ 0$, for all $k$; also, $\hat{q}^k_{i,i}-q_i \rightarrow 0$, for all $i\in \mc{I}$. Intuitively, Theorem~\ref{th:main2} is based on the fact that $\hat{\mc{F}}^k$ approaches $\mc{F}$ in \eqref{eq:opF} asymptotically (i.e., when $\tilde{Q}^k\approx \0$), hence a contractivity property similar to Lemma~\ref{lem:contractivityF} can be ensured for any big-enough $k$. 
Specifically, we note that $({{\Qbar}}+\tilde{Q}^k)^{-1}={{\Qbar}}^{-1}-({\Qbar}({\Qbar}+\tilde{Q}^k))^{-1}\tilde{Q}^k=:{{\Qbar}}^{-1}-P^k$, since the matrices involved are diagonal. Therefore 
$
 \hat{\mc{F}}^k(\x)={\mc{F}}(\x)+\alpha \tilde{\mc{F}}^k(\W \x),
$
with $\mc{F}$ as in \eqref{eq:opF} and $\tilde{\mc{F}}^k(\x):=\mc{R}^\top P^k\F(\x)$. Analogously to Lemma~\ref{lem:Qlipmon}, it can be shown that $\tilde{\mc{F}}^k$ is $\tildeell^k$-Lipschitz in $\H_{\Q}$, with $\tildeell^k:=\eigmax{P^k} \ell \textstyle \sqrt{\eigmax{{Q}}/{\eigmin{{Q}}}}$.
Then, by Lemma~\ref{lem:contractivityF}, $\hat{\mc{F}}$ is $(\sqrt{\rho_\alpha}+\alpha \tildeell^k)$-restricted Lipschitz in $\H_{\Q}$ with respect to $\bs{E}$ (cf. Lemma~\ref{lem:contractivityF}).  Then, analogously to Theorem~\ref{th:main1}, it holds, for all $\k$,  that
\[
\| \x^{k+1}-\x^* \|_{\Q} \leq (\sqrt{\rho_\alpha}+\alpha \tilde{\ell}^k) \| \x^{k}-\x^* \|_{\Q}.
\]
We remark that $\tildeell^k \rightarrow 0$, since 
$\tilde{Q}^k \rightarrow \0$. Hence, for any $\epsilon>0$, the conclusion  follows with $K=(\prod_{k=1}^{\bar{k}} \max\{\sqrt{\rho_\alpha}+\alpha \tilde{\ell}^k,1\})(\sqrt{\rho_\alpha}+\epsilon)^{-\bar{k}}$, where $\bar{k}:=\max\{k \mid \alpha \tilde{\ell}^k>\epsilon\}$. 
 \hfill  \QED 
 
 \smallskip
\subsubsection{Proof of Theorem~\ref{th:main3}}\label{app:th:main3}
	Analogously to the proof of Theorem~\ref{th:main2}, for all $\k$, it holds that $\| \x^{k+1}-\x^*\|_{\Q} \leq \delta^k\| \x^{k}-\x^*\|_{\Q}$, $ \delta^k:= (\sqrt{\rho_{\alpha^k}}+\alpha^k\tilde{\ell}^k)$, with $\rho_{\alpha^k}$ as in \eqref{eq:C1} and $(\tilde{\ell}^k)_{\k}$ a vanishing nonnegative sequence. The conclusion follows because $\prod_{k=0}^{\infty} \delta^k=0$, as we show next.  
	By explicit computation of the quantity in \eqref{eq:C1} and Taylor expansion at $\alpha=0$, it holds,  in a neighborhood
$V_0$ of $\alpha=0$,  that $\sqrt{\rho_\alpha}=1-\bar{\mu}\eigmin{Q}\alpha+o(\alpha)$,
where $o(\alpha)$ is a series of monomial terms at least quadratic in $\alpha$. Take $\bar{k}$ such that, for all $k\geq \bar{k}$, $\tildeell^k\leq \ell^*<\barmu \eigmin{Q}$ for some $\ell^*$, $\alpha^k\in \V_0$ and $\delta^k<1$ (which is always possible, because $\tildeell^k\rightarrow 0$, $\alpha^k\rightarrow 0$ and $\delta^k=1-(\bar{\mu}\eigmin{Q}-\tildeell^k)\alpha^k+o(\alpha^k)$ if $\alpha^k\in\V_0$). Then, $\prod_{k=\bar{k}}^{\infty} \delta^k=0$ if and only if $\sum_{k=\bar{k}}^{\infty} -\log(\delta^k)=\infty$. In turn, by the asymptotic comparison theorem and by the Taylor expansion at $\alpha^k=0$, the latter series diverges if the series $\sum^{\infty}_{k=\bar{k}} \alpha^k (\textstyle {\barmu}\eigmin{Q}-\tildeell^*)$ diverges, which holds by the assumption on $(\alpha^k)_{\k}$. \hfill \QED

\bibliographystyle{IEEEtran}
\bibliography{library}

\begin{thebibliography}{10}
\providecommand{\url}[1]{#1}
\csname url@samestyle\endcsname
\providecommand{\newblock}{\relax}
\providecommand{\bibinfo}[2]{#2}
\providecommand{\BIBentrySTDinterwordspacing}{\spaceskip=0pt\relax}
\providecommand{\BIBentryALTinterwordstretchfactor}{4}
\providecommand{\BIBentryALTinterwordspacing}{\spaceskip=\fontdimen2\font plus
\BIBentryALTinterwordstretchfactor\fontdimen3\font minus
  \fontdimen4\font\relax}
\providecommand{\BIBforeignlanguage}[2]{{%
\expandafter\ifx\csname l@#1\endcsname\relax
\typeout{** WARNING: IEEEtran.bst: No hyphenation pattern has been}%
\typeout{** loaded for the language `#1'. Using the pattern for}%
\typeout{** the default language instead.}%
\else
\language=\csname l@#1\endcsname
\fi
#2}}
\providecommand{\BIBdecl}{\relax}
\BIBdecl

\bibitem{Barrera2015}
J.~{Barrera} and A.~{Garcia}, ``Dynamic incentives for congestion control,''
  \emph{IEEE Transactions on Automatic Control}, vol.~60, no.~2, pp. 299--310,
  2015.

\bibitem{Saad2012}
W.~{Saad}, Z.~{Han}, H.~V. {Poor}, and T.~{Basar}, ``Game-theoretic methods for
  the smart grid: An overview of microgrid systems, demand-side management, and
  smart grid communications,'' \emph{IEEE Signal Processing Magazine}, vol.~29,
  no.~5, pp. 86--105, 2012.

\bibitem{Li_Chen_Dahleh_2015}
N.~{Li}, L.~{Chen}, and M.~A. {Dahleh}, ``Demand response using linear supply
  function bidding,'' \emph{IEEE Transactions on Smart Grid}, vol.~6, no.~4,
  pp. 1827--1838, 2015.

\bibitem{Ghaderi_2014}
J.~Ghaderi and R.~Srikant, ``Opinion dynamics in social networks with stubborn
  agents: Equilibrium and convergence rate,'' \emph{Automatica}, vol.~50,
  no.~12, pp. 3209 -- 3215, 2014.

\bibitem{Shamma_Arslan_2005}
J.~S. {Shamma} and G.~{Arslan}, ``Dynamic fictitious play, dynamic gradient
  play, and distributed convergence to {N}ash equilibria,'' \emph{IEEE
  Transactions on Automatic Control}, vol.~50, no.~3, pp. 312--327, 2005.

\bibitem{DePersisGrammaticoTAC2020}
C.~{De Persis} and S.~{Grammatico}, ``Continuous-time integral dynamics for a
  class of aggregative games with coupling constraints,'' \emph{IEEE
  Transactions on Automatic Control}, vol.~65, no.~5, pp. 2171--2176, 2020.

\bibitem{BelgioiosoGrammatico_ECC_2018}
G.~Belgioioso and S.~Grammatico, ``Projected-gradient algorithms for
  generalized equilibrium seeking in aggregative games are preconditioned
  forward-backward methods,'' in \emph{2018 European Control Conference}, 2018,
  pp. 2188--2193.

\bibitem{SwensonKarXavier_FictitiousPlay_2015}
B.~{Swenson}, S.~{Kar}, and J.~{Xavier}, ``Empirical centroid fictitious play:
  An approach for distributed learning in multi-agent games,'' \emph{IEEE
  Transactions on Signal Processing}, vol.~63, no.~15, pp. 3888--3901, 2015.

\bibitem{Koshal_Nedic_Shanbag_2016}
J.~Koshal, A.~{Nedić}, and U.~V. Shanbhag, ``Distributed algorithms for
  aggregative games on graphs,'' \emph{Operations Research}, vol.~64, pp.
  680--704, 2016.

\bibitem{TatarenkoShiNedic_CDC2018}
T.~{Tatarenko}, W.~{Shi}, and A.~{Nedić}, ``Accelerated gradient play
  algorithm for distributed {N}ash equilibrium seeking,'' in \emph{2018 IEEE
  Conference on Decision and Control (CDC)}, 2018, pp. 3561--3566.

\bibitem{DePersisGrammatico2018}
C.~{De Persis} and S.~Grammatico, ``Distributed averaging integral {N}ash
  equilibrium seeking on networks,'' \emph{Automatica}, vol. 110, p. 108548,
  2019.

\bibitem{YeHu2017}
M.~{Ye} and G.~{Hu}, ``Distributed {N}ash equilibrium seeking by a consensus
  based approach,'' \emph{IEEE Transactions on Automatic Control}, vol.~62,
  no.~9, pp. 4811--4818, 2017.

\bibitem{BelgioiosoNedicGrammatico2020}
G.~{Belgioioso}, A.~{Nedić}, and S.~{Grammatico}, ``Distributed generalized
  {N}ash equilibrium seeking in aggregative games on time-varying networks,''
  \emph{IEEE Transactions on Automatic Control, DOI:10.1109/TAC.2020.3005922},
  2020.

\bibitem{Pavel2018}
L.~{Pavel}, ``Distributed {GNE} seeking under partial-decision information over
  networks via a doubly-augmented operator splitting approach,'' \emph{IEEE
  Transactions on Automatic Control}, vol.~65, no.~4, pp. 1584--1597, 2020.

\bibitem{Bianchi_arXiv_TAC20_PPP}
M.~Bianchi, G.~Belgioioso, and S.~Grammatico, ``Fast generalized {N}ash
  equilibrium seeking under partial-decision information,'' \emph{arXiv
  preprint 2003.09335}, 2020.

\bibitem{DengNian2019}
Z.~{Deng} and X.~{Nian}, ``Distributed generalized {N}ash equilibrium seeking
  algorithm design for aggregative games over weight-balanced digraphs,''
  \emph{IEEE Transactions on Neural Networks and Learning Systems}, vol.~30,
  no.~3, pp. 695--706, 2019.

\bibitem{Bianchi_LCSS2020}
M.~{Bianchi} and S.~{Grammatico}, ``Fully distributed {N}ash equilibrium
  seeking over time-varying communication networks with linear convergence
  rate,'' \emph{IEEE Control Systems Letters}, vol.~5, no.~2, pp. 499--504,
  2021.

\bibitem{SalehisadaghianiPavel2017_nondoublystochastic}
F.~Salehisadaghiani and L.~Pavel, ``Nash equilibrium seeking with non-doubly
  stochastic communication weight matrix,'' \emph{EAI Endorsed Trans. on
  Collaborative Computing}, vol.~4, no.~13, pp. 1-- 11, 2019.

\bibitem{XiMai_DistributedOptimizationonRowStochastic}
C.~{Xi}, V.~S. {Mai}, R.~{Xin}, E.~H. {Abed}, and U.~A. {Khan}, ``Linear
  convergence in optimization over directed graphs with row-stochastic
  matrices,'' \emph{IEEE Transactions on Automatic Control}, vol.~63, no.~10,
  pp. 3558--3565, 2018.

\bibitem{FacchineiPang2007}
F.~Facchinei and J.~Pang, \emph{Finite-dimensional variational inequalities and
  complementarity problems}.\hskip 1em plus 0.5em minus 0.4em\relax Springer,
  2007.

\bibitem{Cenedese2019AsynchronousAT}
C.~{Cenedese}, G.~{Belgioioso}, Y.~{Kawano}, S.~{Grammatico}, and M.~{Cao},
  ``Asynchronous and time-varying proximal type dynamics in multi-agent network
  games,'' \emph{IEEE Transactions on Automatic Control,
  DOI:10.1109/TAC.2020.3011916}, 2020.

\bibitem{Bianchi_ECC20_ctGNE}
M.~{Bianchi} and S.~{Grammatico}, ``A continuous-time distributed generalized
  {N}ash equilibrium seeking algorithm over networks for double-integrator
  agents,'' in \emph{2020 European Control Conference}, 2020, pp. 1474--1479.

\bibitem{Bauschke2017}
H.~H. Bauschke and P.~L. {Combettes}, \emph{Convex analysis and monotone
  operator theory in Hilbert spaces}.\hskip 1em plus 0.5em minus 0.4em\relax
  Springer, 2017, vol. 2011.

\end{thebibliography}
\end{document}